\documentclass[a4paper]{amsart}
\usepackage{amsmath}

\def\n{n}
\def\d{{n+1}}
\def\Torus{\mathbb{T}}

\def\me{\mathbf{g}}
\def\contr{\neg}

\def\ka{K\"{a}hler~}

\def\contr{\rightharpoonup}

\def\rnum{\mathbb{R}}
\def\cnum{\mathbb{C}}
\def\znum{\mathbb{Z}}

\def\ka{K\"{a}hler\;}

\def\omuno{{\omega_1}}
\def\omdue{{\omega_2}}
\def\omd{{\omega_D}}

\def\x12{\mathbf{X_1}\times_B \mathbf{X_2}}
\def\mex12{\me_{_{\x12}}}

\def\me{\mathbf{g}}

\newtheorem{teo}{Theorem}[section]

\newtheorem{cor}[teo]{Corollary}

\newtheorem{lem}[teo]{Lemma}
\newtheorem{pro}[teo]{Proposition}
\newtheorem{dfn}[teo]{Definition}

\newtheorem{rmk}[teo]{Remark}

\title{Self-dual manifolds and mirror symmetry for the quintic threefold}
\author{Michele Grassi}
\date{January 3, 2003}
\begin{document}
\begin{abstract}
For all $m\geq 1$ we build a two-dimensional family of smooth manifolds of real dimension $3m + 2$ and use it  to interpolate between the anticanonical family in $\cnum\mathbb{P}^{m+1}$ and its mirror dual. The main tool is the notion of self-dual manifold.
\end{abstract}
\maketitle
\section{Introduction}
In the present paper we describe a way to interpolate geometrically between the large \ka structure limit point in the (\ka) moduli space of the anticanonical divisor in $\cnum\mathbb{P}^\n$ and a large complex structure limit point in the complex structure moduli space of  its mirror partner (which is a submanifold of the complex manifold $H^\n$). 
The interpolation is achieved by constructing a two dimensional family of smooth manifolds of (real) dimension $3(\n -1) +2$. For instance, for the quintic threefold we obtain a two dimensional family of $11$ dimensional smooth manifolds. These manifolds are endowed with a structure, which we introduced in ~\cite{G2}, and we call a {\em weakly self-dual structure} (or WSD structure for brevity). The definition is given at the beginning of the next section and  involves a Riemannian metric and three smooth $2$-forms. The manifolds depend on two parameters $\rho_1,\rho_2$. Qualitatively, what happens is that fixing $\rho_2$ determines the "shape" of the limiting manifold, while if we let $\rho_1$ go to $+\infty$ we get the large \ka structure limit, and if we let $\rho_1$ go to zero we get the large complex structure limit. Moreover as $\rho_2$ goes to infinity, the limiting manifolds approach in a normalized Gromov-Hausdorff sense the anticanonical divisors of $\cnum\mathbb{P}^\n$ and their mirror duals. Another (dual) construction relates in the same way a  large complex structure limit on the anticanonical divisors  of $\cnum\mathbb{P}^\n$ and the large \ka structure limit of its mirror dual family . To clarify what happens on the boundary of the deformation space, it is useful to imagine the deformation space as a square, with the four sides associated to the values $\rho_1 = +\infty,\rho_2 = \rho_2^{min},\rho_1 = 0,\rho_2 = +\infty$ respectively. Call the first three sides $A,T,B,S$ respectively, and call also $M_A$ the vertex common to the sides $A$ and $S$, and similarly call $M_B$ the vertex common to the sides $B$ and $S$. Then the point $M_A$ corresponds to the large \ka structure limit point in the (\ka) moduli space of the anticanonical divisor in $\cnum\mathbb{P}^\n$ and $M_B$ corresponds to large complex structure limit point in the complex structure moduli space of  its mirror partner. Moreover, the points on the (interior of the) boundary $A$ are "infinitely inflated" $\Torus^\n$ fibrations over the sphere $S^{\n -1}$, the points on the (interior of the) boundary $T$ are complex tori of (complex) dimension $\n$ with a choice of a \ka structure on them, and the points on the (interior of the) boundary $B$ are "infinitely inflated" real $\Torus^\n$'s. Finally, there seems to be no easy interpretation for the interior of the boundary $S$, as the objects that one obtains are wildly singular from a metric point of view.  The distance used to take the limits in the above discussion is normalized Gromov-Hausdorff distance. The above picture of the deformation space for our self-dual manifolds has a striking similarity with the conjectural picture of the moduli space of superconformal field theories described by Kontsevich and Soibelman in ~\cite{KS}. This agreement is in accordance with a more general conjectural picture, in which (weakly) self-dual manifolds can be used to build superconformal field theories via a process similar to a sigma-model construction. However, such a procedure has not yet been established in a mathematically rigorous way even for the more classical Calabi-Yau manifolds. We cannot therefore claim that our construction verifies in any way the conjecures of ~\cite{KS} for the anticanonical families in projective spaces. \\
We should point out that the limits $M_A$ and $M_B$ are not bona fide limits, but more like "asymptotic" limits. Indeed, at $M_A$ for   any choice of large $\rho_2$ we must choose $\frac{e^{2\pi^2\rho_2^2}}{\rho_1}$ small enough  to have that the self-dual manifold converges in  normalized Gromov-Haurdorff distance to  the set of points $[z_0,...,z_\n]\in \cnum\mathbb{P}^\n_{\rho_1}$ which satisfy the equation $\prod_iz_1 = 0$,
where we indicate with $\cnum\mathbb{P}^\n_{\rho_1}$ projective space endowed with the symplectic form which is $\rho_1^2$ times the Fubini-Study one. Similarly at $M_B$  for   any choice of large $\rho_2$ we must have $\rho_1\rho_2$ small enough to have that the self-dual manifold converges in  Gromov-Haurdorff distance to  the set of points $[z_0,...,z_\n]\in H^\n_{\rho_1,\rho_2}$ which satisfy the equation $\prod_iz_1 = 0$
where we indicate with $H^\n_{\rho_1,\rho_2}$ the space $H^\n$ (cf. Definition ~\ref{def-hn})  endowed with a complex structure (defined in Definition ~\ref{dfn:cpxhn}) compatible with  the (induced) Fubini-Study two-form and which for $\rho_1$ which goes to $0$ tends to a "large complex structure limit". As mentioned before, we can also make a dual construction, which interpolates a large compex structure limit point for the anticanonical divisor in projective space with the large \ka structure limit point of its mirror dual in $H^\n$. 

The interpolating manifolds are constructed via a procedure which has a toric flavor to it, and starts from the reflexive polytope associated to $\cnum\mathbb{P}^\n$. The toric nature of the construction is reflected in the fact that the resulting manifolds have a (free) action by the real torus $\Torus^{\n}\times\Torus^{\n}$.  Moreover, the limiting procedure involves a geometric deformation (reflected in a rescaling of the parameter $\rho_1$) which implies the rescaling of the metric on one of the two fibrations by a factor, and on the other fibration by the inverse of the same factor. As mentioned before, depending on  the fact that we let the parameter $\rho_1$ go to zero or to infinity, we approach one or the other limit point of the deformation space. This description of mirror symmetry has some similarity with the conjectural description of the mirror involution contained in the paper ~\cite{SYZ} by Strominger, Yau and Zaslow, although in a (possibly) unexpected way. Indeed, we do not build special lagrangian fibrations on the Calabi-Yau manifolds themselves near the limit points, but we end up with "special" tori fibrations on (higher dimensional) WSD manifolds, which approximate the Calabi-Yau ones only in Gromov-Hausdorff sense. The idea that this could be a way to avoid the complications associated with building special lagrangian fibrations in the geometric approach to mirror symmetry is what led us to the definition of self-dual manifolds in the first place. 
As for $T$-duality (cf. ~\cite{SYZ} for the definition), it does not hold in the manifolds that we build, except possibly in an approximate way near the boundary of the deformation space. However, we think that there should be a way to identify inside the whole deformation space of the manifolds that we build a subspace made up of "$T$-dual" WSD manifolds, for which $T$-duality holds (exactly) for the two $\Torus^\n$ fibrations mentioned above. This would verify the conjecture of ~\cite{SYZ} without necessarily implying the existence of special lagrangian fibrations on the limiting Calabi-Yau manifolds. The most natural way to impose this duality condition on the structure would imply asking for the codimension $\n - 1$ differential form giving the Riemannian volume on the distribution associated to the $\Torus^\n\times\Torus^n$ fibration to be closed. This condition however would determine a differential equation on the components of the metric very similar to the Monge-Amp\`{e}re equation associated to the Calabi-Yau condition, and hence its integration might not be completely trivial to perform.

In the paper ~\cite{G2} we proved that self-dual manifolds can be used to build an interpolating family for mirror pairs of elliptic curves and of Affine-\ka manifolds. In both these cases however the dual special lagrangian fibrations do exist on the Calabi-Yau manifolds, and we actually used these fibrations to build the interpolating self-dual manifolds as fibre products over the common base of the fibrations.  Notice also that in that paper we built self-dual manifolds, which are WSD but enjoy also one more property. Finally, we should mention that our construction generalizes to more general polytopes. It is however not clear what (if any) relevance the resulting WSD manifolds have to mirror symmetry.
We now give a description of the content of the various sections.\\
In section ~\ref{sec:prelim} we very briefly introduce self-dual manifolds and weakly self-dual manifolds. We do not explore any of their properties, as that has been already done in ~\cite{G2}. We need however to generalize slightly the definition that was given in ~\cite{G2}, as in the present paper we need to consider degenerate weakly self-dual  structures, while in the cited reference we only considered the non degenerate case. We then introduce the polytopes associate to projective space when they are considered as toric varieties, and their dual polytopes. We do not state or prove any facts on reflexive polytopes, as the only property that we need to perform the construction is easily stated and proved directly in our case. We feref to ~\cite{B} for more on reflexive polytopes, and for the construction of the conjectural mirror partners using them.\\
In section ~\ref{sec:maincost} we perform the main construction of weakly self-dual manifolds $\mathbb{X}^{\n-1}_{k_1,k_2}$ starting from the polytope $\Delta_\n$ for $\cnum\mathbb{P}^\n$ and its dual $\Delta_\n^*$. We then prove that what we obtained is actually weakly self-dual.\\
In section ~\ref{sec:gpd} we first consider a natural action by the torus $\Torus^\n\times\Torus^\n$ on the manifold $\mathbb{X}^{\n-1}_{k_1,k_2}$, which makes it "toric". For reasons of space we do not try to define what is a toric weakly self-dual manifold, even if all the ingredients for a natural generalization of the standard definition would be present. We then define the two fundamental projection maps $\pi_1$ and $\pi_2$ of a weakly self-dual manifold, in the special case of the manifolds $\mathbb{X}^{\n-1}_{k_1,k_2}$. The map $\pi_1$ takes values in $\cnum\mathbb{P}^\n$, while the map $\pi_2$ takes values in the manifold $H^\n$ associated to the polytope $\Delta_\n^*$. We give equations for the images of both $\pi_1$ and $\pi_2$, which will be needed later. We finally introduce a natural geometric deformation which is present whenever one has a nondegenerate weakly self-dual manifold, and can be induced on the $\mathbb{X}^{\n-1}_{k_1,k_2}$ (which are degenerate) in a natural way. This deformation will play a  crucial r\^ole in the following.\\
In section ~\ref{sec:boundary} we study the boundary of the deformation space of the manifolds 
$\mathbb{X}^{\n-1}_{k_1,k_2}$. To do that we define a normalization of the Gromov-Hausdorff distance, which is useful in the following as we need to compare manifolds with divergent diameter, and we we are only interested in their "shape". Using this distance,  we show among other things that there are two special points on the boundary of the deformation space. One of them corresponds to families converging to the large \ka structure limit of the anticanonical divisor in projective space, while the other corresponds to families converging to a large complex structure limit on its mirror.\\
In the final section we briefly sketch how one can generalize the construction to more general polytopes, and we conclude with some remarks and some questions.

This paper was written while at the Mathematical Sciences Research Institute (M.S.R.I.) in Berkeley, California. I would like to thank the organizers who made my stay there possible, and the staff at the Institute for creating a very pleasant and stimulating environment for doing research.
\section{Preliminary facts on self-dual manifolds and reflexive polytopes}
\label{sec:prelim}
In this section, after some preliminary remarks on self-dual and weakly self-dual manifolds, we prove some elementary facts on reflexive polytopes that we will need in the sequel.
In the present paper we will need a slight generalization of the notions presented in 
~\cite{G2}, so we give
 here briefly the basic definitions, without comments or examples. 
We refer to ~\cite{G2} for those, an for a more extensive and 
detailed introduction to them.
\begin{dfn}
A {\em weakly self-dual manifold} (WSD manifold for brevity) is given by a smooth manifold $X$, together with two smooth $2$-forms $\omega_1,\omega_2$ a Riemannian metric and a third smooth  $2$-form $\omega_D$ (the {\em dualizing} form) on it, which satisfy the following conditions:\\
1) $d\omega_1 = d\omega_2 = d\omega_D = 0$ and the distribution $\omega_1^0+\omega_2^0$ is integrable.\\
2)
For all $p\in X$ here exist an orthogonal basis $dx_1,..,dx_m,dy^1_1,...,dy^1_m,dy^2_1,...,dy^2_m$, $dz_1,...,dz_c,dw_1,...,dw_c$ of  $T_p^*X$ such that the $dx_1,..,dx_m,dy^1_1,...,dy^1_m,dy^2_1,...,dy^2_m$ are orthonormal and at $p$
\[\omega_1 = \sum_{i=1}^m dx_i\wedge dy^1_i,~ ~\omega_2 = \sum_{i=1}^m dx_i\wedge dy^2_i,~ ~\omega_D = \sum_{i=1}^m dy^1_i\wedge dy^2_i + \sum_{i=1}^cdz_i\wedge dw_i\]
Any orthogonal basis of $T_pX$ dual to a basis  of $1$- forms  as above  is said to be {\em adapted} to the structure, or {\em standard}. The number $m$ is the {\em rank} of the structure.
\end{dfn}
For a more intrinsic definition of WSD manifolds the reader should refer to  ~\cite{G2}. Here we have chosen the quickest way to introduce them.
\begin{rmk}
The form $\omega_D$ is symplectic once restricted to $\omega_1^0 + \omega_2^0$. We have therefore that $\omega_D^{dim(X) - m}\not=0$.
\end{rmk}
\begin{dfn}
1) A WSD manifold is {\em nondegenerate} if $dim(\omega_1^0\cap\omega_2^0)_p = 0$ at all points (equivalently if its dimension is $3$ times the rank).\\
2) A  WSD manifold is {\em self-dual} (SD manifold for brevity) if all the leaves of the distribution $\omega_1^0+\omega_2^0$ have volume one (with respect to the volume form induced by the metric)
\end{dfn}
The difference with respect to the definitions given in ~\cite{G2} is that there we only considered the nondegenerate case, where $\omega_1^0\cap\omega_2^0 = (0)$. The present definitions of self-dual and  weakly self-dual manifold simplify  to those ones in this special case. In the following we will be mainly interested in the case where $dim(\omega_1^0\cap\omega_2^0)_p = 2$ at all $p$, and hence $dim(X) = 3m +2$.
In the nondegenerate case, condition $2$ in the definition is enough to determine $\omega_D$ starting from $\omega_1,\omega_2$ and the metric. It is however not true that in this case the properties of $\omega_1,\omega_2$ and $g$ are enough to guarantee that $\omega_D$ is closed.  The above definitions are all that we will  need from ~\cite{G2}.\\
Let us now come to reflexive polytopes. For their definition and their basic properties we refer to  ~\cite{B}.
\begin{dfn}
Let $\Delta_\n$ be the polytope associated to $\mathbb{P}^\n$, when considered as a toric variety in the standard way. Indicate with $\Delta_\n^*$ the dual polytope
\end{dfn}
The polyhedra $\Delta_\n, \Delta_n^*$ are  given by \\
$\Delta_\n = \text{convex hull of}\left\{v_1 = (\n,-1,....,-1),.., v_\n = (-1,,...,n), v_{\n+1} = (-1,...,-1)\right\}$\\
$\Delta_n^* = \text{convex hull of}\left\{u_1 = (1,....,0),.., u_\n = (0,...,1), u_{\n+1} = (-1,...,-1)\right\}$\\
Point  1) and the corresponding part of point 3) of the following definition are standard constructions, and can be found for example in ~\cite{Gu}.
\begin{dfn} In the notations of the previous definition, define:\\
1)  The linear maps of real vector spaces
 $F_{\Delta_\n}, F_{\Delta_\n^*}:\rnum^{\d}\to \rnum^\n$ are
 \[F_{\Delta_\n}(x_1,...,x_\d) = \sum_{i=1}^\d x_i u_i,~~F_{\Delta^*_\n}(x_1,...,x_\d) = \sum_{i=1}^\d x_i v_i\]
2) The maps $F_{\Delta_\n}^*, F_{\Delta_\n^*}^*: \rnum^\n \to \rnum^d$ are the maps obtained from the transposes of the matrices of  of $F_{\Delta_\n}, F_{\Delta_\n^*}$ with respect to the standard bases.\\
 3) The group morphism  $f_{\Delta_\n},~f_{\Delta_\n^*}~:~\Torus^\d\to\Torus^\n$ and $f_{\Delta_\n}^*,~f_{\Delta_\n^*}^*~:~\Torus^\n\to\Torus^\d$ are the maps induced by the 
$F_{\Delta_\n}, F_{\Delta_\n^*}, F_{\Delta_\n}^*, F_{\Delta_\n^*}^*$ respectively, after quotienting by the integer lattices of the spaces.
 \end{dfn}  
 Notice that in the definition of $F_{\Delta_\n}$ you use the vertices of $\Delta_\n^*$, and viceversa.
The maps can also be defined explicitely using the standard bases as
\[
\begin{array}{l}
F_{\Delta_\n}((x_1,...,x_{\d}) )=  (x_1-x_{\d},...,x_\n-x_{\d})\\
F_{\Delta_\n^*}((x_1,...,x_{\d}) )=  (\n x_1-\sum_{i\not= 1}x_i,...,\n x_\n-\sum_{i\not=\n}x_i)
\end{array}\]
\begin{lem}
1) The linear map $F_{\Delta_\n}F_{\Delta_\n^*}^*  = F_{\Delta_\n^*}F_{\Delta_\n}^*$  is $\d$ times the  identity of $\rnum^\n$\\
2) $Ker\left(F_{\Delta_\n}\right) 
\oplus Im\left(F_{\Delta_\n^*}^*\right) = 
Ker\left(F_{\Delta_\n^*}\right)
\oplus \left(F_{\Delta_\n}^*\right) = \rnum^\d$
\end{lem}
\textit{Proof}
1) The proof  is an easy direct computation.\\
2) This follows from part $1$ and the fact that the rank of all the maps $F_{\Delta_\n}$, $F_{\Delta_\n^*}^*$, $F_{\Delta_\n^*}$, $F_{\Delta_\n}^*$ is $\n$.
\qed

The first point of the following definition is a standard object, described for example in ~\cite{Gu}.
\begin{dfn}
1) We indicate with $N_{\Delta_\n}$ (resp. $N_{\Delta_\n^*}$) the kernel of $f_{\Delta_\n}$ (respectively of $f_{\Delta_\n^*}$). \\
2) Define  
\[D_{\Delta_\n} = N_{\Delta_\n}\cap Im (f_{\Delta_\n^*}^*),~ ~D_{\Delta_\n^*} = N_{\Delta_\n^*} \cap Im (f_{\Delta_\n}^*)\]
\end{dfn}
\begin{cor}
\label{kertorus}
\[\begin{array}{ll}
1) & Ker\left(f_{\Delta_\n}f_{\Delta_\n^*}^*\right) = Ker\left(f_{\Delta_\n^*}f_{\Delta_\n}^*\right) = \left(\znum/(\d)\right)^\n \subset \Torus^\n\\
2) & N_{\Delta_\n} + Im\left(f_{\Delta_\n^*}^*\right) = N_{\Delta_\n^*} +  \left(f_{\Delta_\n}^*\right) = \Torus^\d
\end{array}\]
\end{cor}

\begin{rmk}
$D_{\Delta_\n} = f_{\Delta_\n^*}^*\left(Ker\left(f_{\Delta_\n}f_{\Delta_\n^*}^*\right)\right)$ and 
$D_{\Delta_\n^*} = f_{\Delta_\n}^*\left(Ker\left(f_{\Delta_\n^*}f_{\Delta_\n}^*\right)\right)$, therefore they are both finite groups.
\end{rmk}
\begin{lem} 
The maps $f_{\Delta_\n^*}^*, f_{\Delta_\n}$ (resp. $f_{\Delta_\n}^*, f_{\Delta_\n^*}$) induce group isomorphisms
\[\Torus^\d/N_{\Delta_\n} \cong f_{\Delta_\n^*}^*\left(\Torus^\n\right)/D_{\Delta_\n} \cong \Torus^\n /\left(\znum/(\d)\right)^\n,\]
\[\Torus^\d/N_{\Delta_\n^*} \cong f_{\Delta_\n}^*\left(\Torus^\n\right)/D_{\Delta_\n^*} \cong \Torus^\n /\left(\znum/(\d)\right)^\n\]
\end{lem}
\textit{Proof} This is clear from the definition of the $D_{\Delta_\n}, D_{\Delta_\n^*}$ and the previous corollary.
\qed

The reader familiar with Delzant's construction of toric varieties via a symplectic reduction might be a little surprised by the previous statement. To clarify it, one should notice that the above isomorphisms translate into isomorphisms of the fibres of the reduction {\em away from the points where the original fibres collapse}. In this respect, the existence of the isomorphism above is linked to  the fact that there is an open subset inside any toric variety which is isomorphic to a product of $\cnum^*$'s.\\
\section{From the polytopes $\Delta_\n$ to weakly self-dual manifolds}
\label{sec:maincost}
In this section we give the main construction, which starts from the reflexive polytopes $\Delta_\n$ and $\Delta_\n^*$ and builds a (two dimensional family of) compact smooth WSD manifolds of dimension $3(\n -1) + 2$. There is also a dual construction (which we don't write down explicitely), in which the r\^oles of $\Delta_\n$ and of $\Delta_\n^*$ get switched.\\
We indicate with $(r_0,...,r_\n,\theta_0,...,\theta_\n) = (\bar{r},\bar{\theta})$ the standard coordinates on $(\cnum^*)^\d$, so that the standard holomorphic coordinates are $x_j + i y_j = z_j = r_j e^{2\pi i \theta_j}$ for $j\in\{0,...,\n\}$. Then the map $\mu(r_0,...,r_\n,\theta_0,...,\theta_\n) = -(\pi r_0^2,...,\pi r_\n^2)$ is  a moment map for the standard $\Torus^\d$ action on $(\cnum^*)^\d$ (given by translation of the $\theta_j$'s in the previous notation). We are using the symplectic form $\sum_{j = 0}^\n 2\pi r_j dr_j\wedge d\theta_j = \frac{-1}{2 i}\sum_{j = 0}^\n dz_j\wedge d\bar{z}_j$ which is compatible with the standard flat metric of $\cnum^* \cong \rnum^2\setminus \{0\}$, given in coordinates by $\sum_j \left((dr_j)^2 + 4\pi^2 r_j^2 (d\theta_j)^2\right)$.
\begin{dfn}
Let $\mu:\cnum^\d\to \rnum^\d$ be the moment map described above. Let $(\cnum^\d)^*\times_\mu(\cnum^\d)^*$ be the associated fibred product space. By using two copies 
$(\bar{r},\bar{\theta})$ and $(\bar{r},\bar{\eta})$ of the standard coordinates on $(\cnum^*)^\d$, 
we obtain global coodinates $(\bar{\theta},\bar{r},\bar{\eta})$ on $(\cnum^\d)^*\times_\mu(\cnum^\d)^*$. We put on it the self-dual structure 
$((\cnum^\d)^*\times_{\rnum^\d} (\cnum^\d)^*,\tilde{\omega}_1,\tilde{\omega}_2,\tilde{\mathbf{g}},\tilde{\omega}_D)$
given in coordinates by 
\[\tilde{\omega}_1 = 2\pi\sum_{i=0}^{\n}r_i dr_i\wedge d\theta_i,~\tilde{\omega}_2 = \frac{1}{2\pi}\sum_{i=0}^{\n}\frac{1}{r_i} dr_i\wedge d\eta_i,\]
\[\tilde{\mathbf{g}} = \sum_{i=1}^{\d}\left(dr_i^2 + 4\pi^2r_i^2d\theta_i^2 + \frac{1}{4\pi^2 r_i^2}d\eta_i^2\right),~ ~
\tilde{\omega}_D = \sum_{i=0}^\n d\theta_i \wedge d\eta_i\]
\end{dfn}
\begin{lem}
The above structure on $(\cnum^\d)^*\times_\mu(\cnum^\d)^*$ is self-dual
\end{lem}
\textit{Proof}
The two-forms involved are clearly smooth and closed, and the distribution $\tilde{\omega}_1^0 + \tilde{\omega}_2^0$ is the span of the $\frac{\partial}{\partial \theta_i},\frac{\partial}{\partial \eta_j}$ and is therefore clearly integrable. The basis
\[\frac{1}{2\pi r_0}\frac{\partial}{\partial \theta_0},...,\frac{1}{2\pi r_\n}\frac{\partial}{\partial \theta_\n},
\frac{\partial}{\partial r_0},...,\frac{\partial}{\partial r_\n},  
2\pi r_0\frac{\partial}{\partial \eta_0},...,2\pi r_\n\frac{\partial}{\partial \eta_\n}\]
is by inspection orthonormal and adapted to the structure. This proves that the structure is WSD. To prove that it is actually self-dual, we observe that the coordinates $(\bar{\theta},\bar{\eta})$ provide an identification of the leaves of the distribution $\tilde{\omega}_1^0 + \tilde{\omega}_2^0$ with the torus $\Torus^\d\times\Torus^\d$. Under this identification, for fixed $\bar{r}$, the metric gets sent to a flat metric on $\Torus^\d\times\Torus^\d$, which gives length $2\pi r_0,...,2\pi r_\n,\frac{1}{2\pi r_0},...,\frac{1}{2\pi r_0}$ to the $2(\d)$ $\Torus^1$-factors of $\Torus^{2(\d)}$. From this, the volume of such a leaf is 
$\left(\prod_{i=0}^\n (2\pi r_i)\right)\left(\prod_{j=0}^\n (\frac{1}{2\pi r_j})\right) = 1$.
\qed

The WSD manifolds that we are going to build are "morally" (poly)symplectic reductions of $(\cnum^\d)^*\times_\mu(\cnum^\d)^*$ by the action defined below. From this point of view, the map  $(\mu_1,\mu_2)$ of the next definition plays the r\^ole of the moment map. However, the presence of a section simplifies things in our context, and allows us to build them as quotients by a finite group of  submanifolds of $(\cnum^\d)^*\times_\mu(\cnum^\d)^*$, without having to develop the theory of reduction for such actions. Such a theory is  in our opinion interesting, but describing it would take us too far from our present objective.
\begin{dfn}
1) The group $N_{\Delta_\n}\times N_{\Delta_\n^*}$ (and hence also its subgroup $D_{\Delta_\n}\times D_{\Delta_\n^*}$) acts naturally and freely on $(\cnum^\d)^*\times_\mu(\cnum^\d)^*$, via its inclusion in $\Torus^\d\times\Torus^\d$.  Moreover, this action has a section, given by the points $\tilde{\sigma} = \left\{(\bar{\theta},\bar{r},\bar{\eta})~|~\bar{\theta} = \bar{\eta} = 0\right\}$.\\
2) 
The map $(\mu_1,\mu_2)~:~(\cnum^\d)^*\times_\mu(\cnum^\d)^*\to\rnum^2$ is defined  as
\[(\mu_1,\mu_2)(\bar{\theta},\bar{r},\bar{\eta}) = -\left(\pi\sum_{i=0}^\n r_i^2,\frac{1}{2\pi}log\left(\prod_{i=0}^\n r_i\right)\right) \]
\end{dfn}
\begin{lem}
1) The action of  $N_{\Delta_\n}\times N_{\Delta_\n^*}$  preserves  the $2$-forms  and the metric of the self-dual structure of $(\cnum^\d)^*\times_\mu(\cnum^\d)^*$, and the two maps $\mu_1,\mu_2$ are equivariant with respect to it. \\
2) $\mu_1^{-1}(k_1)\cap \mu_2^{-1}(k_2)\subset \cnum^\d\times_\mu\cnum^\d$ is non-empty if and only if $\frac{-k_1}{\pi}e^{\frac{4\pi}{\d}k_2} \geq \d$\\
3) For $ \frac{-k_1}{\pi}e^{\frac{4\pi}{\d}k_2} > \d $,  for any point $p\in\mu_1^{-1}(k_1)\cap \mu_2^{-1}(k_2)\subset \cnum^\d\times_\mu\cnum^\d$ we have $d\mu_1\wedge d\mu_2\not= 0$
\end{lem}
\textit{Proof}
1) Clear by inspection.\\
2) The minimum of $\mu_2$ for a fixed value of $\mu_1$ is obtained when all the $r_i$ are equal to a common value $r$. At such a point we have 
\[k_2 = \mu_2 = -\frac{\d}{2\pi}log(r),~ ~k_1 = \mu_1 = -\pi(\d)r^2 = -\pi(\d)e^{-\frac{4\pi}{\d}k_2}\]
and hence $\frac{-k_1}{\pi}e^{\frac{4\pi}{\d}k_2} = \d$. This proves that whenever $\frac{-k_1}{\pi}e^{\frac{4\pi}{\d}k_2} \geq \d$ the two equations $\mu_1 = k_1, \mu_2 = k_2$ admit a common solution, while  when $\frac{-k_1}{\pi}e^{\frac{4\pi}{\d}k_2} < \d$ there cannot be any point where they are both satisfied.\\
3) Assume that for some $p\in\mu_1^{-1}(k_1)\cap \mu_2^{-1}(k_2)$ we have $d\mu_1\wedge d\mu_2 = 0$. Then there exists $t\in\rnum^+$ such that 
\[\forall i ~ ~2\pi r_i = \frac{t}{2\pi r_i}\]
Therefore $\forall i ~ ~r_i = \frac{\sqrt{t}}{2\pi}$, and 
\[k_1 = -\frac{(\d)t}{4\pi^2},~ ~k_2 = -\frac{\d}{4\pi}log\left(\frac{t}{4\pi^2}\right)\] 
From this it follows that $\frac{-k_1}{\pi}e^{\frac{4\pi}{\d}k_2} = \d$, whatever value $t$ had.
\qed
\begin{dfn}
Let
\[\rho_1 = \sqrt{\frac{-k_1}{\pi}},~ ~
\rho_2 = \sqrt{\frac{\d}{4\pi^2}log\left(\frac{-k_1}{\pi}\right) + \frac{k_2}{\pi}}\]
\end{dfn}
The quantities $\rho_1,\rho_2$ will show up again many times. We will among other things use them as the fundamental parameters to describe the deformations of the manifolds that we are about to build. In terms of $\rho_2$, the condition above becomes 
$e^{\frac{4\pi^2 \rho_2^2}{\d}} \geq \d$.
\begin{dfn}
Assume that $e^{\frac{4\pi^2\rho_2^2}{\d}} > \d$. Then the space $\mathbb{X}_{k_1,k_2}^{\n -1}$ is defined as a quotient of a submanifold  of $\tilde{\mathbb{X}}_{k_1,k_2}^{\n -1}\subset (\cnum^\d)^*\times_\mu(\cnum^\d)^*$:
\[\tilde{\mathbb{X}}_{k_1,k_2}^{\n -1}~=~\left\{ (\bar{\theta},\bar{r},\bar{\eta}) \in \left(Im(f_{\Delta_\n^*}^*)\times Im(f_{\Delta_\n}^*)\right)\sigma(\bar{r})~|~
\mu_1(\bar{\theta},\bar{r},\bar{\eta}) = k_1,\mu_2(\bar{\theta},\bar{r},\bar{\eta}) = k_2\right\}\]
and 
\[\mathbb{X}_{k_1,k_2}^{\n -1}~=~\tilde{\mathbb{X}}_{k_1,k_2}^{\n -1}/D_{\Delta_\n}\times D_{\Delta_\n^*}\]
\end{dfn}
Notice that the action by the group $D_{\Delta_\n}\times D_{\Delta_\n^*}$ is fixed-point free, so the quotient is a smooth (compact) manifold as soon as $(k_1,k_2)$ is a regular value for $(\mu_1,\mu_2)$.
\begin{lem}
The submanifold $\tilde{\mathbb{X}}_{k_1,k_2}^{\n -1}$ is an integral manifold for the distribution  of vectors annichilated by the space of forms
\[<~\left(\sum_{i=0}^\n\frac{\partial}{\partial\theta_i}\right)\contr \tilde{\omega}_D,~ \left(\sum_{i=0}^\n\frac{\partial}{\partial\eta_i}\right)\contr \tilde{\omega}_D,~d\mu_1,~d\mu_2~>\]
\end{lem}
\textit{Proof}
A tangent vector is of the form
$\sum_i a_i \frac{\partial}{\partial r_i} + df_{\Delta_\n^*}^*(u) + df_{\Delta_\n}^*(w)$. We have only to verify that $df_{\Delta_\n^*}^*(u)\in \left(\sum_{i=0}^\n\frac{\partial}{\partial\theta_i}\contr \tilde{\omega}_D\right)^0$, $df_{\Delta_\n}^*(w)\in \left(\sum_{i=0}^\n\frac{\partial}{\partial\eta_i}\contr \tilde{\omega}_D\right)^0$, as the other vanishings are clear. However, by construction $\sum_{i=0}^\n\frac{\partial}{\partial\theta_i}\in Ker(df_{\Delta_\n^*})$ and hence its dual is in the orthogonal to the image of the transpose of the differential with respect to any metric. Here we are just asking this for the metric determined by the basis $\frac{\partial}{\partial\theta_i}$ (recall that contrancting with $\omega_D$ is equivalent to sending to the dual as the form is orthonormal symplectic with respect to the metric induced by the bases). Notice that we are identifying the orbits of the two copies of $\Torus^d$ without mentioning it. The same argument proves the other inclusion, and we have therefore the thesis.
\qed
\begin{teo}
In the previous notations, assume  that 
$\frac{-k_1}{\pi}e^{\frac{4\pi}{\d}k_2} > \d$, so that at  any point 
$p\in\mu_1^{-1}(k_1)\cap \mu_2^{-1}(k_2)\subset \cnum^\d\times_\mu\cnum^\d$ we have $d\mu_1\wedge d\mu_2\not= 0$. Then the self-dual structure on 
$\cnum^\d\times_\mu\cnum^\d$ induces three $2$-forms 
$\omega_1,\omega_2,\omega_D$ and a metric  
$\mathbf{g}$ on the tangent space at $p$ for any 
$p\in\mathbb{X}_{k_1,k_2}$, plus a section $\sigma$. The forms and the metric determine  a WSD structure on $\mathbb{X}_{k_1,k_2}^{\n -1}$
\end{teo}

\textit{Proof} It is clearly enough to prove that the self-dual structure on the ambient space induces a WSD structure on $\tilde{\mathbb{X}}_{k_1,k_2}^{\n -1}$, as the quotient by a finite group of elements which preserve the structure preserves the WSD property. Let 
\[X_1 = \sum_{i=0}^\n\frac{\partial}{\partial\theta_i},~ ~
X_2 = \frac{1}{4\pi^2}\sum_{i=0}^\n\frac{1}{r_i^2}\frac{\partial}{\partial\theta_i},~ ~
Y_1 = \sum_{i=0}^\n\frac{\partial}{\partial\eta_i},~ ~
Y_2 = 4\pi^2\sum_{i=0}^\n r_i^2\frac{\partial}{\partial\eta_i}\]
Then 
$X_1\contr\tilde{\omega}_1 = Y_2\contr\tilde{\omega}_2$ and $X_2\contr\tilde{\omega}_1 = Y_1\contr\tilde{\omega}_2$, 
$\| X_1 ^2\| = \| Y_2\|^2 = 4\pi^2\sum_{i=0}^\n r_i^2$ and 
$\| X_2 \|^2= \| Y_1\|^2 = \frac{1}{4\pi^2}\sum_{i=0}^\n \frac{1}{r_i^2}$. 
We also have $<X_1, X_2> = <Y_1,Y_2> = \d$.

We have, indicating with $X^*$ the $1$-form dual to the vector $X$ with respect to the metric, that 
\[X_1^* = \sum_{i=0}^\n 4\pi^2r_i^2 d\theta_i,~ ~
X_2^* = \sum_{i=0}^\n d\theta_i,~ ~
Y_1^* = \sum_{i=0}^\n\frac{1}{4\pi^2r_i^2} d\eta_i,~ ~
Y_2^* = \sum_{i=0}^\n d\eta_i\]
By construction,
\[\left(Y_1\contr\tilde{\omega_D}\right)^0 = (X_2^*)^0 = <X_2>^\perp,~ ~\left(X_1\contr\tilde{\omega_D}\right)^0 = (Y_2^*)^0 =  <Y_2>^\perp\]
and therefore from the previous lemma the tangent space at any point to $\tilde{\mathbb{X}}_{k_1,k_2}^{\n -1}$ is just $T_p\tilde{\mathbb{X}}_{k_1,k_2}^{\n -1} = \left(<X_2,Y_2>^\perp ~\cap~ Ker(d\mu_1)~\cap ~Ker(d\mu_2)\right)_p$. We need to show that the forms of the structure restricted to this space satisfy the pointwise conditions for a WSD structure (part $2$ of the definition).\\
As $d\mu_{1}\wedge d\mu_{2}\not= 0$, we have that the dimension of $Ker\left(d\mu_{1}\right)\cap Ker\left(d\mu_{2}\right)$ is $3\d-2$. Furthermore, for $i\in\{1,2\}$
\[\tilde{\omega}_1^0\oplus\tilde{\omega}_2^0\subset Ker\left(d\mu_{1}\right)\cap Ker\left(d\mu_{2}\right),
~ ~dim\left(Ker\left(d\mu_{i}\right)\cap
\left(\tilde{\omega}_1^0\oplus\tilde{\omega}_2^0\right)^\perp\right) = \n\]
It follows that the intersection $Ker\left(d\mu_{1}\right)\cap Ker\left(d\mu_{2}\right)\cap \left(\tilde{\omega}_1^0\oplus\tilde{\omega}_2^0\right)^\perp$ has dimension $\n -1$.
Take an orthonormal basis $v_1,...,v_{\n-1}$ of this space, and complete it to orthonormal bases (for $i \in\{1,2\}$)
\[<v_1,...,v_{\n-1},v_n^i> = Ker\left(d\mu_{i}\right)\cap
\left(\tilde{\omega}_1^0\oplus\tilde{\omega}_2^0\right)^\perp,\]
\[<v_1,...,v_{\n-1},v_n^i,v_{\n+1}^i> = \left(\tilde{\omega}_1^0\oplus\tilde{\omega}_2^0\right)^\perp\]
For $i\in  \{1,2\}$, let  $u_1^i,...,u_{\n+1}^i$ (resp. $w_1^i,...,w_{\n+1}^i$) be the  basis inside $\tilde{\omega}_2^0\subset (\tilde{\omega}_1^0)^\perp$ (resp. inside $\tilde{\omega}_1^0\subset (\tilde{\omega}_2^0)^\perp$) dual  to the basis $v_1,...,v_{\n-1},v_\n^i,v_{\n +1}^i$ with respect to the two-form $\tilde{\omega}_1$ (resp. $\tilde{\omega}_2$), which is symplectic when rectricted to the space $(\tilde{\omega}_1^0)^\perp$ (resp. the space $(\tilde{\omega}_2^0)^\perp$ ). Then by construction 
\[u_\d^1 = \frac{1}{\|X_1\|}X_1,~ ~ u_\d^2 = \frac{1}{\|X_2\|}X_2,~ ~w_\d^2 = \frac{1}{\|Y_1\|}Y_1,~ ~ w_\d^1 = \frac{1}{\|Y_2\|}Y_2 \]
It follows that at all points $p$ the orthogonal set of (nonzero) vectors
\[X_1 - \frac{<X_1,X_2>}{\|X_2\|^2}X_2,u_1^1,...,u_{\n-1}^1,v_1,...,v_{\n -1},w_1^2,...,w_{\n-1}^2, Y_1 - \frac{<Y_1,Y_2>}{\|Y_2\|^2}Y_2\]
spans the space orthogonal to $X_2,Y_2$ and is inside $\left(Ker(d\mu_1)~\cap ~Ker(d\mu_2)\right)_p$, and therefore is a basis for the tangent space to $\tilde{\mathbb{X}}_{k_1,k_2}^{\n -1}$ at the point $p$. Notice that the quantity $<X_1,X_2>^2 - \| X_2\|^2\| X_1\|^2$ (and therefore both the vectors $X_1 - \frac{<X_1,X_2>}{\|X_2\|^2}X_2$ and $Y_1 - \frac{<Y_1,Y_2>}{\|Y_2\|^2}Y_2$) is always different from zero, as for it to vanish the vectors $X_1,X_2$ would have to be collinear, and that never happens as long as $d\mu_1\wedge d\mu_2 \not= 0$. 
The subset $u_1^1,...,u_{\n-1}^1,v_1,...,v_{\n -1},w_1^2,...,w_{\n-1}^2$ is an orthonormal set and
\[\omega_1|_{<X_2,Y_2>^\perp} = \sum_{i=1}^{\n -1} (v_i)^*\wedge (u_i^1)^*,~ ~\omega_2|_{<X_2,Y_2>^\perp} = \sum_{i=1}^{\n -1} (v_i)^*\wedge (w_i^2)^*\]
This guarantees that the forms $\omuno,\omdue$ and the metric satisfy the axioms required for a WSD structure regarding the pointwise conditions. For the form $\omd$, first observe that from the definition of $\tilde{\omega}_D$ we get 
\[\omd|_{<u_1^1,...,u_{\n-1}^1,v_1,...,v_{\n -1},w_1^2,...,w_{\n-1}^2>} = \sum_{i=1}^{\n -1} (u_i^1)^*\wedge (w_i^2)^*\]
To prove that $\omd^\n \not= 0$ it is therefore enough to show that it is not zero when restricted to the span of $X_1 - \frac{<X_1,X_2>}{\|X_2\|^2}X_2, Y_1 - \frac{<Y_1,Y_2>}{\|Y_2\|^2}Y_2$. As the vectors $X_1,X_2,Y_1,Y_2$ form a basis for the space that they generate, there must be coefficients $\{a_{ij}\}$ such that
\[{\tilde{\omega}_D}|_{<X_1,X_2,Y_1,Y_2>} = a_{11}X_1^*\wedge Y_1^* + a_{12}X_1^*\wedge Y_2^* + a_{21}X_2^*\wedge Y_1^* + a_{22}X_2^*\wedge Y_2^*\]
To determine the $\{a_{ij}\}$ we use the defining conditions for $X_2$ and $Y_2$ and the fact that the structure on $(\cnum^\d)^*\times_\mu(\cnum^\d)^*$ is nondegenerate (and hence $\tilde{\omega}_D$ is uniquely determined):
\[Y_2^* = \tilde{\omega}_D(X_1, -) =\| X_1\|^2 a_{11} Y_1^* + \| X_1\|^2 a_{12} Y_2^* + (\n +1)a_{21} Y_1^* + (\n +1)a_{22} Y_2^*\]
\[Y_1^* = \tilde{\omega}_D(X_2, -) = (\n+1)a_{11} Y_1^* + (\n+1)a_{12} Y_2^* + \| X_2\|^2 a_{21} Y_1^* + \| X_2\|^2a_{22} Y_2^*\]
From this we get the system
\[\left\{\begin{array}{l}
\| X_1\|^2 a_{11} + (\n +1)a_{21}  = 0 \\
\| X_1\|^2 a_{12} + (\n +1)a_{22} = 1\\
(\n+1)a_{11} +  \| X_2\|^2 a_{21}  = 1\\
(\n+1)a_{12} +  \| X_2\|^2a_{22} = 0
\end{array}\right.\]
which can be successively reduced to 
\[\left\{\begin{array}{l}
a_{21}  = -\frac{\| X_1\|^2}{(\n +1)} a_{11} \\
\left((\n +1)^2  - \| X_1\|^2 \| X_2\|^2\right)a_{22}  = (\n +1) \\
\left((\n +1)^2 - \| X_2\|^2\| X_1\|^2 \right)a_{11}  =  (\n +1)\\
a_{12} = -\frac{ \| X_2\|^2}{(\n+1)}a_{22}
\end{array}\right.
\iff
\left\{\begin{array}{l}
a_{11} = \frac{\n + 1}{(\n +1)^2 - \| X_2\|^2\| X_1\|^2}\\
a_{12} = -\frac{\| X_2\|^2}{(\n +1)^2 - \| X_2\|^2\| X_1\|^2}\\ 
a_{21}  = -\frac{\| X_1\|^2}{(\n +1)^2 - \| X_2\|^2\| X_1\|^2}\\
a_{22} = \frac{\n + 1}{(\n +1)^2 - \| X_2\|^2\| X_1\|^2}

\end{array}\right.\]
Then $\left(<X_1,X_2>^2 - \| X_2\|^2\| X_1\|^2\right){\tilde{\omega}_D}|_{<X_1,X_2,Y_1,Y_2>}$ has the expression
\[(\d)X_1^*\wedge Y_1^* - \| X_2\|^2X_1^*\wedge Y_2^* - \| X_1\|^2X_2^*\wedge Y_1^* + (\d)X_2^*\wedge Y_2^*\]
and hence 
\[\omd\left(X_1 - \frac{\d}{\|X_2\|^2}X_2, Y_1 - \frac{\d}{\|Y_2\|^2}Y_2\right) = \tilde{\omega}_D\left(X_1 - \frac{\d}{\|X_2\|^2}X_2, Y_1 - \frac{\d}{\|X_1\|^2}Y_2\right) = \] 
\[ \frac{\d}{(\d)^2 - \|X_1\|^2\| X_2\|^2}\left(1 - \frac{\| X_2\|^2}{\|X_1\|^2}  - \frac{\| X_1\|^2}{\|X_2\|^2} +  
\frac{(\d)^2}{\|X_1\|^2\| X_2\|^2}\right) = \frac{\d}{\|X_1\|^2\| X_2\|^2}\]
Therefore the given basis of $<X_2,Y_2>^\perp$ is adapted to the WSD structure.
The forms $\omuno,\omdue,\omd$ are also closed by construction . The remaining  condition to verify is that the distribution $\omuno^0 + \omdue^0$ is integrable. This follows from the fact that it is the restriction to $\tilde{\mathbb{X}}_{k_1,k_2}^{\n -1}$ of an integrable distribution in the ambient space.
 \qed
 \begin{rmk}
 From the proof of the previous theorem, and computing the norms of the vectors $X_1 - \frac{\d}{\|X_2\|^2}X_2$ and $Y_1 - \frac{\d}{\|Y_2\|^2}Y_2$, we obtain that the norm of the form $\omd$ when restricted to the span of the degenerate distribution is $\frac{\d}{\|X_1\|^2\| X_2\|^2 - (\d)^2}$. 
 \end{rmk}
 \begin{rmk}
 The forms $X_1^*$ and $Y_1^*$ induce on $\mathbb{X}^{\n -1}_{k_1,k_2}$ a structure which is a generalization of a (nondegenerate) contact structure in the same way as the forms $\omega_1,\omega_2$ are a generalization of a (degenerate) symplectic structure. Using this analogy as a guide, one could use this "polycontact" structure to induce a nondegenerate WSD structure on $\mathbb{X}^{\n -1}_{k_1,k_2}\times \Torus^2$. 
\end{rmk}
\begin{dfn}
We indicate with $\mathbb{X}^{\n-1}_{k_1,k_2}$ both the quotient manifold constructed above and the same manifold endowed  with the WSD structure $(\omuno,\omdue,\me,\omd)$. When referring to a generic element of the family, we may drop the subscripts $k_1,k_2$.
\end{dfn}
\begin{rmk}
We have $dim(\mathbb{X}^{m}) = 3m +2$. In particular $dim(\mathbb{X}^{3}) = 11$
\end{rmk}
\section{Natural group actions, projections and deformations}
\label{sec:gpd}
The following action is the same natural one that one obtains on a toric variety after building it via symplectic reduction starting from its (dual) polytope. In our case,  we chose to use a section for the group action by which one reduces, instead of performing the quotient, which makes the proof different from the standard one.
\begin{lem}
There is a natural  free action of the group $\Torus^\n\times\Torus^\n$ on $\mathbb{X}^{\n -1}$, 
induced by the morphism $f_{\Delta_\n}\times f_{\Delta^*_\n}~:~\Torus^\d\times \Torus^\d\to \Torus^\n\times\Torus^\n$. The orbits of the first factor are the leaves of the foliation $\omega_1^0$, while the orbits of the second factor are the leaves of the foliation $\omega_1^0$.
\end{lem}
\textit{Proof}
Take $(t_1,t_2)\in \Torus^\n\times\Torus^\n$, and let $p = [\bar{\theta},\bar{r},\bar{\eta}]\in \mathbb{X}^{\n -1}$. To induce the action, we use $(s_1,s_2)\in \Torus^\d\times \Torus^\d$ such that $f_{\Delta_\n}\times f_{\Delta^*_\n}(s_1,s_2) = (t_1,t_2)$. Pick also $\tilde{p} = (\bar{\theta},\bar{r},\bar{\eta}) \in\tilde{\mathbb{X}}^{\n -1}$ which  maps to $p$ under the natural projection from $\tilde{\mathbb{X}}^{\n -1}$ to $\mathbb{X}^{\n -1}$, $p = [\tilde{p}]$. As $f_{\Delta_\n^*}^*(\Torus^\n) + N_{\Delta_\n} = \Torus^\d$ and $f_{\Delta_\n}^*(\Torus^\n) + N_{\Delta_\n^*} = \Torus^\d$ (from Corollary ~\ref{kertorus}), there must be an 
$u\in N_{\Delta_\n}\times N_{\Delta_\n^*}$ such that $u(s_1,s_2)\tilde{p} \in f_{\Delta_\n^*}^*(\Torus^\n)\times f_{\Delta_\n}^*(\Torus^\n)\tilde{\sigma} $. Define 
\[(t_1,t_2)p = [u(s_1,s_2)\tilde{p}]\]
We have to check that this is a well posed definition (and then it is automatically a group action), and that the action so defined is free. For the first, observe that the ambiguity in the choice of $\tilde{p}$ is associated to the possible multiplication by an element 
$(u_1,u_2)\in D_{\Delta_\n}\times D_{\Delta_\n^*}$, while the ambiguity in the choice of $(s_1,s_2)$ is associated to the possible multiplication by an arbitrary element 
$(v_1,v_2) \in N_{\Delta_\n}\times N_{\Delta_n^*}$. Now if $u(s_1,s_2)\tilde{p} \in f_{\Delta_\n^*}^*(\Torus^\n)\times f_{\Delta_\n}^*(\Torus^\n)\tilde{\sigma}$, we have that $(u(u_1^{-1}v_1^{-1},u_2^{-1}v_2^{-1})(u_1v_1,u_2v_2)(s_1,s_2)\tilde{p} \in f_{\Delta_\n^*}^*(\Torus^\n)\times f_{\Delta_\n}^*(\Torus^\n)\tilde{\sigma}$. Summing up, all the choices made are summed up in the ambiguity in the choice of the element $u$. If however also  $v(s_1,s_2)\tilde{p} \in f_{\Delta_\n^*}^*(\Torus^\n)\times f_{\Delta_\n}^*(\Torus^\n)\tilde{\sigma}$ for another element $v\in N_{\Delta_\n}\times N_{\Delta_\n^*}$, we have that 
\[uv^{-1} \in \left( f_{\Delta_\n^*}^*(\Torus^\n)\times f_{\Delta_\n}^*(\Torus^\n)\right) \cap \left(N_{\Delta_\n}\times N_{\Delta_n^*}\right) = D_{\Delta_\n}\times D_{\Delta_\n^*}\]
This implies that $[u\tilde{p}] = [v\tilde{p}]$, and we have a well defined action. To see that the action is free, assume in the previous notation that 
$[u(s_1,s_2)\tilde{p}] = p$. This implies that there must be $(a_1,a_2)\in D_{\Delta_\n}\times D_{\Delta_\n^*}$ with $(a_1,a_2)u(s_1,s_2)\tilde{p} = \tilde{p}$. As the original $\Torus^\d\times\Torus^\d$ action is free, this implies that $(a_1,a_2)u(s_1,s_2) = (e,e)$ (the identity element). This however implies that $(s_1,s_2)\in N_{\Delta_\n}\times N_{\Delta_n^*}$ and therefore $(t_1,t_2) = (e,e)$ as required.\\
The last statement is clear by inspection.
\qed
\begin{teo}
\label{pione}
There is a smooth map $\pi_1$ from  
$\mathbb{X}^{\n -1}_{k_1,k_2}$ to $\cnum\mathbb{P}^\n$, with the leaves of the distribution $\omega_1^0$ as fibres. This map is equivariant with respect to the natural 
$\Torus^\n\times \Torus^\n$ action on 
$\mathbb{X}^{\n -1}_{k_1,k_2}$ and the natural $\Torus^\n$ action on 
$\cnum\mathbb{P}^\n$, and with respect to the projection onto the first factor
$pr_1: \Torus^\n\times \Torus^\n \to \Torus^\n$. Moreover, the form $\omega_1$ induces via this map a \ka form  on 
$\cnum\mathbb{P}^\n$ which is $\rho_1^2$ times the Fubini-Study one. The image of $\pi_1$ is the set of points $[z_0,...,z_\n]\in \cnum\mathbb{P}^\n$ which satisfy the equation 
\[\prod_{i=0}^\n(z_i\bar{z}_i) = \frac{1}{e^{4\pi^2\rho_2^2}}\left(\sum_{i=0}^\n z_i\bar{z}_i\right)^{\d}\]
\end{teo}
\textit{Proof}
We use the identification (as a symplectic manifold with a torus action)
\[\cnum\mathbb{P}^\n_\lambda = \left\{(z_0,...,z_n)\in (\cnum^\d)^*~|~\sum_{i=o}^\n \|z_i\|^2 = \lambda\right\}/N_{\Delta_\n}\]
provided by Delzant's construction. As $\lambda$ varies in $\rnum^+$ we get different symplectic forms $\omega_{(\lambda)}$ on the space, and we have $\omega_{(\lambda)} = \lambda\omega_{(1)}$. Moreover, a direct computation shows that $\omega_{(1)}$ is just the Fubini-Study \ka form on $\cnum\mathbb{P}^\n$. There is then a natural way of defining $\pi_1$, namely if $\lambda = \frac{-k_1}{\pi}$, 
\[\pi_1([\bar{\theta},\bar{r},\bar{\eta}]) = [r_0e^{2\pi i\theta_0},...,r_\n e^{2\pi i\theta_\n}]\]
First of all we have to verify that this definition is well posed, so choose 
$(u_1,u_2)$ in  $D_{\Delta_\n}\times D_{\Delta_\n^*}$. We then have $ [r_0e^{2\pi i\theta_0},...,r_\n e^{2\pi i\theta_\n}] =  [r_0e^{2\pi iu_1(\theta_0)},...,r_\n e^{2\pi iu_1(\theta_\n)}]$ because $D_{\Delta_\n} \subset N_{\Delta_\n}$ by construction. So the map $\pi_1$ is well defined. If $\pi_1([\bar{\theta},\bar{r},\bar{\eta}]) = \pi_1([\bar{\theta}^\prime,\bar{r}^\prime,\bar{\eta}^\prime])$, then $\bar{\theta} = u\bar{\theta}^\prime$ for some $u\in N_{\Delta_\n}$. As however both $\bar{\theta}$ and $\bar{\theta}^\prime$ lie in the same lateral class of $\Torus^\d/D_{\Delta_\n}$ with respect to $N_{\Delta_\n}/D_{\Delta_\n}$, it must be $u\in D_{\Delta_\n}$ and therefore $[\bar{\theta},\bar{r},\bar{\eta}] = [\bar{\theta}^\prime,\bar{r}^\prime,\bar{\eta}]$. This implies that the fibres of $\pi_1$ are obtained varying $\bar{\eta}$, and therefore they coincide with the integral manifolds of $\omega_1^0$.\\
Let $(t_1,t_2)\in \Torus^\n$. From the last statement in the previous lemma, we have to verify that $\pi_1(t_1[\bar{\theta},\bar{r},\bar{\eta}]) = t_1[r_0e^{2\pi i\theta_0},...,r_\n e^{2\pi i\theta_\n}]$. Recall that the action of $t_1$ on $\cnum\mathbb{P}^\n$ is obtained by first lifting it to $s_1\in\Torus^\d$ along $f_{\Delta_\n}$, and then applying $s_1$ to any lifting of the point in $\cnum\mathbb{P}^\n$ to a point in $(\cnum^\d)^*$ (cf. ~\cite{Gu} for details). Comparing this with the definition of the action on $\mathbb{X}^{\n -1}_{k_1,k_2}$, the statement is clear by inspection.\\
For the statement concerning the symplectic forms, it is enough to prove that the pull-back of the symplectic form $\omega_{(\lambda)}$ of $\cnum\mathbb{P}^\n_\lambda$ (with $\lambda = \frac{-k_1}{\pi} = \rho_1^2$) along the composed map $\tilde{\pi}_1: \tilde{\mathbb{X}}^{\n -1}_{k_1,k_2}\to \cnum\mathbb{P}^\n_\lambda$ is the form $\omega_1$ on $\tilde{\mathbb{X}}^{\n -1}_{k_1,k_2}$. This however is clear, as both forms are induced (in one case via the inclusion and in the other via a projection) by the symplectic form of $(\cnum^\d)^*$. \\
The last statement is clear by inspection, once we translate the symplectic reduction presentation of projective space into the standard one as a quotient by a  $\cnum^*$ action.
\qed

The following map should be thought of as similar to the Cremona transformation which sends the anticanonical divisor of projective space to itself, having  its irreducible component divisors and their interserctions switched. 
\begin{dfn}
\label{def:defphi}
The map $\phi_{\rho_1,\rho_2}$ from $(\cnum^*)^\d\times_\mu(\cnum^*)^\d $ to itself is defined as
\[\phi_{\rho_1,\rho_2}(\bar{\theta},\bar{r},\bar{\eta}) = (
\bar{\theta},
\rho_1e^{-2\pi^2 \rho_2^2r_0^2},...,
\rho_1e^{-2\pi^2 \rho_2^2r_\n^2},
-\bar{\eta})\] 
\end{dfn}
\begin{lem}
\label{defphi}
The map  $\phi_{\rho_1,\rho_2}$ is a diffeomorphism from $(\cnum^*)^\d\times_\mu(\cnum^*)^\d $ onto the set of points
$\left\{
(\bar{\theta},\bar{r},\bar{\eta})\in 
(\cnum^*)^\d\times_\mu(\cnum^*)^\d
~|~\forall i~r_i < \rho_1
\right\} $. Moreover, 
\[\phi_{\rho_1,\rho_2}^*(\tilde{\omega}_1) = \sum_{i=0}^\n \frac{8\pi^3\rho_1^2\rho_2^2 r_i}{ e^{4\pi^2\rho_2^2 r_i^2}} dr_i\wedge d\theta_i,~ ~  
\phi_{\rho_1,\rho_2}^*(\tilde{\omega}_2) = 2\pi\rho_2^2\sum_{i=0}^\n r_i dr_i\wedge d\eta_i \]
\[\phi_{\rho_1,\rho_2}^*(\tilde{\mathbf{g}}) = 
\sum_{i=0}^\n \frac{4\pi^2 \rho_1^2}{e^{4\pi^2\rho_2^2r_i^2}}(d\theta_i)^2 + 
\sum_{i=0}^\n\frac{16\pi^4\rho_1^2\rho_2^4r_i^2}{e^{4\pi^2\rho_2^2r_i^2}}(dr_i)^2 + 
\frac{e^{4\pi^2\rho_2^2r_i^2}}{4\pi^2 \rho_1^2}(d\eta_i)^2\]
\[\phi_{\rho_1,\rho_2}^*(- k_1 + \mu_1) = (-k_1)\left(1 - \sum_{i=0}^\d e^{-4\pi^2\rho_2^2r_i^2}\right), 
\phi_{\rho_1,\rho_2}^*(\mu_2 - k_2) = \pi\rho_2^2\left(\sum_{i=0}^\n (r_i)^2  - 1\right) \]
The tensor $\tilde{J}_2  = \sum_i\left(2\pi r_i\frac{\partial}{\partial \eta_i}\otimes dr_i - \frac{1}{2\pi r_i}\frac{\partial}{\partial r_i}\otimes d\eta_i\right) $ on $(\cnum^\d)^*\times_\mu(\cnum^\d)^*$ gets pulled back  by the isomorphism $\phi_{\rho_1,\rho_2}$  to the expression 
\[J_2 =  \sum_i\left(\frac{8\pi^3r_i\rho_1^2\rho_2^2}{ e^{4\pi^2\rho_2^2r_i^2}}\frac{\partial}{\partial \eta_i}\otimes dr_i - \frac{e^{4\pi^2\rho_2^2r_i^2}}{8\pi^3r_i\rho_1^2\rho_2^2}\frac{\partial}{\partial r_i}\otimes d\eta_i\right)\]
\end{lem}
\textit{Proof}
The proof of the lemma is just an easy direct computation.
\qed

We recall below the definition of the (singular) manifold associated to the polytope dual to that of projective space. The reader is advised to consult ~\cite{B}, ~\cite{Gu} and ~\cite{CDGP} for further details on it and on its r\^ole in mirror symmetry for anticanonical divisors of projective space.
\begin{dfn}
\label{def-hn} 
The space $H^\n$ is defined as the quotient of $\cnum\mathbb{P}^\n$ by the action of the finite group $\left(\znum/(\n -1)\right)^\d$ . It can also be presented (as a toric \ka manifold) via the symplectic reduction associated (by Delzant's method) to the polytope $\Delta_\n^*$ dual to the polytope associated to $\cnum\mathbb{P}^\n$.
\end{dfn}
The proof of the following theorem is very similar to that of the previous one, with $\Delta_\n$ replaced by $\Delta_\n^*$. The only real difference is that here we first need to use the isomorphism $\phi$ defined in the previous lemma, to put the structure in a more standard form.
\begin{teo}
\label{pitwo}
There is a smooth map $\pi_2$ from  
$\mathbb{X}^{\n -1}_{k_1,k_2}$ to $H^\n$, with the leaves of the distribution $\omega_2^0$ as fibres. This map is equivariant with respect to the natural 
$\Torus^\n\times \Torus^\n$ action on 
$\mathbb{X}^{\n -1}_{k_1,k_2}$ and the natural $\Torus^\n$ action on 
$\cnum\mathbb{P}^\n$, and with respect to the projection onto the second factor
$pr_2: \Torus^\n\times \Torus^\n \to \Torus^\n$. Moreover, the form $\omega_2$ induces via this map a \ka form  on 
$H^\n$ which is $\rho_2$ times the standard one induced on it by the Fubini-Study form on complex projective space. Using the homogeneous coordinates from the covering projective space, the image of $\pi_2$ is the set of points $[z_0,...,z_\n]\in H^\n$ which satisfy the equation 
\[\sum_{i=0}^\d e^{-4\pi^2 z_1\bar{z}_i} = 1\]
\end{teo}
\textit{Proof}
As mentioned in the definition of $H_\n$ preceding the statement of the theorem, we use the identification (as a symplectic manifold with a torus action)
\[H^\n_\lambda = \left\{(z_0,...,z_n)\in (\cnum^\d)^*~|~\sum_{i=o}^\n \|z_i\|^2 = \lambda\right\}/N_{\Delta_\n^*}\]
provided by Delzant's construction. As $\lambda$ varies in $\rnum^+$ we get different symplectic forms $\omega_{(\lambda)}$ on the space, and we have $\omega_{(\lambda)} = \lambda\omega_{(1)}$. Moreover, a direct computation shows that $\omega_{(1)}$ is just the \ka form on $H^\n$ induced by the  Fubini-Study \ka form on $\cnum\mathbb{P}^\n$. 
We define explicitely $\pi_2$, for $\lambda = \rho_2^2$, as
\[\pi_2([\bar{\theta},\bar{r},\bar{\eta}]) = [\sqrt{\frac{1}{2\pi^2\rho_2^2}log\left(\frac{\rho_1}{r_0}\right)}e^{-2\pi i\eta_0},...,\sqrt{\frac{1}{2\pi^2\rho_2^2}log\left(\frac{\rho_1}{r_\n}\right)} e^{-2\pi i\eta_\n}]\]
To prove that the map is well defined, and that it has the properties claimed in the statement, we consider the map $\pi_1\phi_{\rho_1,\rho_2}$ (where $\phi_{\rho_1,\rho_2}$ is the isomorphism defined in Definition ~\ref{def:defphi}), and proceed to consider this composed map. As $\phi$ is an isomorphism on an open set containing the points used in the definition of $\mathbb{X}^{\n -1}_{k_1,k_2}$, this is harmless.  We  indicate with $\omuno^\prime,\tilde{\omega_1}^\prime, \mu_2^\prime$ etc. the pulled back structures along $\phi_{\rho_1,\rho_2}$, to avoid confusing them with the original ones. We then have from the previous lemma that $\tilde{\omega_2}^\prime = 2\pi\rho_2^2\sum_{i=0}^\n r_1 dr_i\wedge d\theta_i$, $\mu_2^\prime -k_2 = \pi\rho_2^2\left(\sum_{i=0}^\n r_i^2 - 1\right)$ and the new composed map is in expressed in coordinates as
\[\pi_2\phi_{\rho_1,\rho_2}([\bar{\theta},\bar{r},\bar{\eta}]) = [r_0e^{2\pi i\eta_0},...,r_\n e^{2\pi i\eta_\n}]\]
At this point it is clear that the proof that the map is well defined and equivariant with respect to the torus actions, that the fibres are the the leaves of the distribution $\omega_2^0$ and that the form $\omega_2$ induces $\rho_2^2$ times the Fubini-Study form on the target space is an exact replica of the proof of the analogous facts contained in the Theorem ~\ref{pione}. We won't reproduce the argument here, as it would mean simply interchanging the indices one and two and the polytope $\Delta_\n$ with $\Delta_\n^*$ everywhere.

It remains to be verified the equation for the image of the map. This however is proved in Lemma ~\ref{defphi}, where it is stated that 
\[\phi^*(- k_1 + \mu_1) = (-k_1)\left(1 - \sum_{i=0}^\d e^{-4\pi^2\rho_2^2r_i^2}\right)\]
\qed

We will need the following remark and proposition in the next section.
\begin{rmk}
The expression
\[\tilde{\mathbf{g}}_{\rho_1,\rho_2} = \sum_{i=0}^\n \frac{4\pi^2 \rho_1^2}{e^{4\pi^2\rho_2^2r_i^2}}(d\theta_i)^2 + 
\sum_{i=0}^\n\frac{16\pi^4\rho_1^2\rho_2^4r_i^2}{e^{4\pi^2\rho_2^2r_i^2}}(dr_i)^2 + 
\frac{e^{4\pi^2\rho_2^2r_i^2}}{4\pi^2 \rho_1^2}(d\eta_i)^2\]
defines a  nondegenerate metric  on $(\cnum^*)^\d\times_\mu(\cnum^*)^\d $ 
\end{rmk}
\qed
\begin{dfn}
\label{dfn:cpxhn}
We indicate with $J_{\lambda_1,\lambda_2}$ the tensor induced on the complement of  $\left\{[\bar{r},\bar{\eta}]~|~\prod_ir_i = 0\right\}$ inside $H^\n$ by the tensor
\[\sum_i\left(\frac{8\pi^3r_i\lambda_1^2\lambda_2^2}{ e^{4\pi^2\lambda_2^2 r_i^2}}\frac{\partial}{\partial \eta_i}\otimes dr_i - \frac{e^{4\pi^2 \lambda_2^2 r_i^2}}{8\pi^3r_i\lambda_1^2\lambda_2^2}\frac{\partial}{\partial r_i}\otimes d\eta_i\right)\]
defined on  $(\cnum^*)^\d$, after quotienting by $N_{\Delta_\n^*}$ the set $\left\{(\bar{r},\bar{\eta})~|~\sum_ir_i^2 = 1\right\}$.
\end{dfn}

The following definition is a generalization of one found in ~\cite{G2}
\begin{dfn}
The deformation $\alpha_t\left(\mathbb{X}^m_{k_1,k_2}\right)$ of $\mathbb{X}^m_{k_1,k_2}$ is defined as the WSD manifold obtained by the same procedure as that used for $\mathbb{X}^m_{k_1,k_2}$, by reducing the self-dual structure
\[(\tilde{\omega}_1)_t = 2\pi t^2\sum_{i=0}^{\n}r_i dr_i\wedge d\theta_i,~(\tilde{\omega}_2)_t = \frac{1}{2\pi}\sum_{i=0}^{\n}\frac{1}{r_i} dr_i\wedge d\eta_i,\]
\[\tilde{\mathbf{g}}_t = \sum_{i=1}^{\d}\left(t^2dr_i^2 + 4\pi^2r_i^2 t^2d\theta_i^2 + \frac{1}{4\pi^2 t^2 r_i^2}d\eta_i^2\right),~ ~(\tilde{\omega}_D)_t = \sum_{i=0}^\n d\theta_i \wedge d\eta_i\]
on $(\cnum^\d)^*\times_\mu(\cnum^\d)^*$ with respect to the same moment maps $\mu_1= -\pi\sum_{i=1}^{\d}r_i^2,~\mu_{2} = \frac{-1}{2\pi}\log\left(\prod_{i=1}^{\d}r_i\right)$ and the same group action as before
\end{dfn}
\begin{pro}
For $t > 0$, 
\[\mathbb{X}^{m}_{t^2k_1,k_2 - \frac{\d}{2\pi}logt}\cong \alpha_t\left(\mathbb{X}^{m}_{k_1,k_2}\right)\]
\end{pro} 
\textit{Proof}
Consider the smooth $\psi_t$  map from  $(\cnum^\d)^*\times_\mu(\cnum^\d)^*$ to itself given by
\[\psi_t(\bar{\theta},\bar{r},\bar{\eta}) = (\bar{\theta},t\bar{r},\bar{\eta})\]
Then
\[\psi_t^*(\tilde{\omega}_1) = (\tilde{\omega}_1)_t,~ ~\psi_t^*(\tilde{\omega}_2) = (\tilde{\omega}_2)_t,~ ~\psi_t^*(\tilde{\mathbf{g}}) = \tilde{\mathbf{g}}_t,~ ~\psi_t^*(\tilde{\omega}_D) = (\tilde{\omega}_D)_t\]
Moreover,
\[\psi^{-1}\left\{(\bar{\theta},\bar{r},\bar{\eta})~|~-\pi\sum_{i=1}^{\d}r_i^2 = t^2k_1,~\frac{-1}{2\pi}\log\left(\prod_{i=1}^{\d}r_i\right) = k_2 - \frac{\d}{2\pi}log(t)\right\} = \]
\[\left\{(\bar{\theta},\bar{r},\bar{\eta})~|~-\pi\sum_{i=1}^{\d}r_i^2 = k_1,~\frac{-1}{2\pi}\log\left(\prod_{i=1}^{\d}r_i\right) = k_2\right\}\]
and $\psi$ is equivariant with respect to the group action, so it induces an isomorphism of the reduced spaces:
\[\bar{\psi}~:~ \alpha_t\left(\mathbb{X}^{m}_{k_1,k_2}\right) \to 
\mathbb{X}^{m}_{t^2k_1,k_2 - \frac{\d}{2\pi}log(t)}\]
\qed
\begin{cor}
If we use the notation $\alpha_t(k_1,k_2) = (t^2k_1, k_2 - \frac{\d}{2\pi}log(t))$, we have that 
\[\rho_1(\alpha_{t}(k_1,k_2)) = t\rho_1,~ ~\rho_2(\alpha_{t}(k_1,k_2)) = \rho_2\]
\end{cor}
\begin{lem}
The correspondence $(k_1,k_2)\to (\rho_1,\rho_2)$ is a smooth bijection on the subset formed by the points where $\mathbb{X}_{k_1,k_2}^{\n -1}$ is well defined. This set corresponds to the set $\{(\rho_1,\rho_2)~|~e^{\frac{4\pi^2\rho_2^2}{\d}} > \d\}$ in the $(\rho_1,\rho_2)$ space
\end{lem}
\textit{Proof}
The subset formed by the points where $\mathbb{X}_{k_1,k_2}^{\n -1}$ is well defined is simply $\{(k_1,k_2)~|~\frac{-k_1}{\pi}e^{\frac{4\pi}{\d}k_2} > \d\}$, which corresponds under the map to the set in the statement of the lemma. The map is also cleary bijective when restricted to this domain and to this codomain, with a smooth inverse.
\qed

\section{The boundary of the deformation space}
\label{sec:boundary}
In this section we analyze the boundary points of the deformation space of $\mathbb{X}^{\n -1}_{k_1,k_2}$, as we vary $k_1$ and $k_2$ (or equivalently $\rho_1$ and $\rho_2$).
To to that in a quantitative way, we have first define  a normalized Gromov-Hausdorff distance which is well suited for our purposes. When reading the definition, recall that our goal is to compare the manifolds  $\mathbb{X}_{k_1,k_2}^m$ with large symplectic structure limit points of symplectic manifolds and large complex structure limit points of complex manifolds.  A similar choice of normalization  was made in ~\cite{GW} to study large complex structure limit points of $K3$ surfaces. We indicate with $diam(X)$ the diameter of the compact metric space $X$, that is $max\{d(x,y)~|~x,y\in X\}$.
\begin{dfn}
Let $d_{GH}$ denote the usual Gromov-Hausdorff distance between metric spaces. \\
The {\em normalized Gromov-Hausdorff distance} $d_{NGH}$ between the compact metric spaces $M,N$ not both with zero diameter is 
\[d_{NGH}(M,N) = \frac{2d_{GH}(M,N)}{diam(M) + diam(N)}\]
We also define $d_{NGH}(pt,pt) = 0$.
\end{dfn}
The intuitive idea is that NGH distance for manifolds with divergent diameter is like "looking at them from a distance". We won't need this here, but notice that NGH distance is interesting also for manifolds with diameter tending to zero. In particular, the space with only one point is isolated in the topology induced by $d_{NGH}$.
\begin{rmk}
For $t\in \rnum^+$ and $M$ a metric space we indicate with $tM$ the metric space with the same underlying topological space, with a distance function wich is obtained by multiplying the distance on $M$ by $t$.  Then for all metric spaces $M,N$ we have that $d_{NGH}(tM,tN) = d_{NGH}(M,N)$
\end{rmk}
The following two theorems state that it is possible to approximate very precisely (in the NGH sense) a family of Calabi-Yau hypersurfaces approaching the large \ka structure limit point (respectively the large complex structure limit point) in their moduli space. They don't say that it is possible to approximate precisely any given Calabi-Yau manifold, however. In this sense, they say that we can  {\em approximate the families}, although we do not attempt to formalize this last concept, as what we mean should be clear from the statement of the theorems.
\begin{teo}
\label{NGH1}
For any $\epsilon > 0$ there exist $R$ such that for all $\rho_2 > R$ and for  $\frac{e^{2\pi^2\rho_2^2}}{\rho_1}$ small enough:\\
1) The fibres of the map $\pi_1$ have diameter bounded by a constant multiple of $\frac{e^{2\pi^2\rho_2^2}}{\rho_1}$.\\
2)  The Riemannian manifold $\mathbb{X}^{\n -1}_{k_1,k_2}$ is at a (standard) Gromov-Hausdorff distance of less than $\epsilon$  from its image under the map $\pi_1$\\
3) The  image under the map $\pi_1$ of $\mathbb{X}^{\n -1}_{k_1,k_2}$ is at a NGH distance less than $\epsilon$ from the hypersurface $\prod_{i=0}^\n z_i = 0$ inside $\cnum\mathbb{P}^\n_{\rho_1}$
\end{teo}
\textit{Proof}
As we will not seek optimal constants for the approximation, we will prove the statements for the space $\tilde{\mathbb{X}}^{\n -1}_{k_1,k_2}$ and the map $\tilde{pi}_1$ instead. This will clearly prove also  the analogous statements for $\mathbb{X}^{\n -1}_{k_1,k_2}$ and $\pi_1$.

1) 
From the definition of $\tilde{\mathbb{X}}^{\n -1}_{k_1,k_2}$, we see that it is formed by points $(\bar{\theta},\bar{r},\bar{\eta})$ which satisfy the equations $\sum_1 r_i^2 = \rho_1^2$, $\prod_i r_i = e^{-2\pi k_2}$.
It follows that for fixed $r_j$, the maximum value which $\prod_{i\not= j} r_i$ can assume is obtained when all $r_i, ~i\not=j$ are equal to $\sqrt{\frac{\rho_1^2 - r_j^2}{\n}}$. Therefore the smallest value which a single $r_j$ can assume is
\[m_r = \frac{e^{-2\pi k_2}}{\left(\sqrt{\frac{\rho_1^2 - r_j^2}{\n}}\right)^\n} >  \frac{\n^{\frac{\n}{2}}}{\rho_1^\n e^{2\pi k_2}} = \frac{\n^{\frac{\n}{2}}\rho_1}{e^{2\pi^2\rho_2^2}}\]
and hence the biggest value which $r_j^{-1}$ can assume is 
$\n^{-\frac{\n}{2}}e^{2\pi^2\rho_2^2}\rho_1^{-1}$. As the fibres of the projection map $\pi_1$ are totally geodetic submanifolds of the torus $\Torus^{\n +1}$ (with a flat invariant metric), and the component tori  $\Torus^1$  of of $\Torus^\d$ are of length $(2\pi r_j)^{-1}$ for varying $j$, from the previous estimate we deduce that the diameter of the fibre is bounded by 
$\pi\n^{-\frac{\n -1}{2}}e^{2\pi^2\rho_2^2}\rho_1^{-1}$. This proves the first point, as the first part is clear by inspection.

2) This is an immediate consequence of the first part. We omit the easy details.

3) We need to show that the set  $\left\{[\bar{\theta},\bar{r}]~|~\prod_i r_i = 0\right\}$ has normalized distance converging to zero to  $\tilde{\pi}_1\left(\tilde{\mathbb{X}}^{\n -1}_{k_1,k_2}\right)$.
For $\delta\in\rnum^+$, define $U_j^\delta, U^\delta \subset \mathbb{X}^{\n -1}_{k_1,k_2}$ as
\[U_j^\delta = \left\{[\bar{\theta},\bar{r}] \in \tilde{\pi}_1\left(\tilde{\mathbb{X}}^{\n -1}_{k_1,k_2}\right) ~|~\forall i\not= j ~r_j > \delta\rho_1 \right\},~ ~U^\delta = \bigcup_j U^\delta_j\]
It is clear that for $\rho_2$ large enough and $\frac{e^{2\pi^2\rho_2^2}}{\rho_1}$ small enough the diameter of both the set above, the set $\left\{[\bar{\theta},\bar{r}]~|~\prod_i r_i = 0\right\}$ and the set $\tilde{\pi}_1\left(\tilde{\mathbb{X}}^{\n -1}_{k_1,k_2}\right)$ is included in the  interval $[\frac{\rho_1}{4},2\pi\rho_1]$.
Then clearly as $\delta$ goes to zero (and for $\rho_2$ large enough and $\frac{e^{2\pi^2\rho_2^2}}{\rho_1}$ small enough) $U^\delta$ converges NGH to  $\tilde{\pi}_1\left(\tilde{\mathbb{X}}^{\n -1}_{k_1,k_2}\right)$. Moreover,
if $[\bar{\theta},\bar{r}] \in U_j^\delta$, then from $\mu_2 = k_2$ we deduce that
\[\frac{r_j^2}{\rho_1^2} = \frac{e^{-4\pi k_2}}{\rho_1^2\prod_{i\not= j}r_i^2} \leq \frac{e^{-4\pi k_2}}{\rho_1^2\delta^\n\left(\frac{-k_1}{\pi}\right)^{\n}} = \frac{1}{\delta^\n e^{4\pi^2\rho_2^2}}\]
This shows that any point in $U^\delta$ us at a distance from a point in 
$\left\{ [\bar{\theta},\bar{r}]~|~\prod_i r_i = 0\right\}$ 
which can me bade to be an arbitrarily small fraction of the diameter by taking $\rho_2$ large enough. Conversely, if $[\bar{\theta},\bar{r}]\in \left\{\prod_i r_i = 0\right\}$, then by rescaling the $r_i$ with $i\not= j$ by a factor smaller thatn one and increasing $r_j$ accordingly to preserve the condition $\sum_i r_i^2 = \rho_1^2$, we can get to a point such that $\mu_2(\bar{r}) = k_2$, which is then in $\tilde{\pi}_1\left(\tilde{\mathbb{X}}^{\n -1}_{k_1,k_2}\right)$. To obtain this we can resize by a factor which is a small fraction of $\rho_1$ (for $\rho_2$ large enough), and therefore we proved that  $\tilde{\pi}_1\left(\tilde{\mathbb{X}}^{\n -1}_{k_1,k_2}\right)$  converges in the NGH distance to
$\left\{[\bar{\theta},\bar{r}]~|~\prod_i r_i = 0\right\}$ for $\rho_2\to\infty$ when $k_2\to \infty$ and $\delta\to 0$. This concludes the proof.
\qed
\begin{pro}
The hypersurface $\prod_{i=0}^\n z_i = 0$ inside $\cnum\mathbb{P}^\n_{\rho_1}$ can be approximated arbitrarily well by a smooth Calabi-Yau hypersurface inside $\cnum\mathbb{P}^\n_{\rho_1}$ with respect to the normalized distance induced by the metric associated to $\rho_1^2$ times the Fubini-Study two-form (equivalently, with respect to the distance induced by the Fubini-Study two-form).
\end{pro}
\label{antic-smooth-1}
\textit{Proof}
We will go quickly over this argument, as its main point is clear and probably in some form or another it is already contained in the literature. The smooth hypersurface can be defined by homogeneous the equation 
\[\prod_{i=0}^\n z_1 = \frac{1}{e^{4\pi^2\rho_2^2}}\sum_{i=0}^\n z_i^{\n +1}\]
which back in non-homogeneous coordinates where $\sum_ir_i^2 = \rho_1^2$ becomes the equation 
\[f_{k_1,k_2}(\bar{\theta},\bar{r}) = \prod_j r_j - e^{-4\pi^2\rho_2^2 - 2\pi i\sum_j\theta_j}\sum_j r_j^{\n+1}e^{2\pi(\n + 1)i\theta_i} = 0\]
For $\delta\in\rnum^+$, define $V_j^\delta, V^\delta \subset \cnum\mathbb{P}^\n_{\rho_1}$ as
\[V_j^\delta = \left\{[\bar{\theta},\bar{r}] \in  \cnum\mathbb{P}^\n_{\rho_1} ~|~f_{k_1,k_2}(\bar{\theta},\bar{r}) = 0,~  ~\forall i\not= j ~r_j > \delta\rho_1 \right\},~ ~V^\delta = \bigcup_j U^\delta_j\]
It is clear that for $\rho_2$ large enough the diameter of both the set above and the sets $\left\{[\bar{\theta},\bar{r}]~|~\prod_i r_i = 0\right\}$, $\left\{[\bar{\theta},\bar{r}] \in  \cnum\mathbb{P}^\n_{\rho_1} ~|~f_{k_1,k_2}(\bar{\theta},\bar{r}) = 0 \right\}$ is $\sim\rho_1$. Further, as $\rho_2$ increases and $\delta$ goes to zero, the   set $V^\delta$ is at a  distance which is an arbitrarily small fraction of $\rho_1$ from  the set $\left\{[\bar{\theta},\bar{r}]~|~\prod_i r_i = 0\right\}$. Moreover it is clear by construction that as $\rho_2$ increases and $\delta$ goes to zero the set $V^\delta$ converges in normalized distance to $\left\{[\bar{\theta},\bar{r}] \in  \cnum\mathbb{P}^\n_{\rho_1} ~|~f_{k_1,k_2}(\bar{\theta},\bar{r}) = 0 \right\}$. This argument shows that the normalized distance of the two sets $\left\{[\bar{\theta},\bar{r}] \in  \cnum\mathbb{P}^\n_{\rho_1} ~|~f_{k_1,k_2}(\bar{\theta},\bar{r}) = 0 \right\}$ and $\left\{[\bar{\theta},\bar{r}]~|~\prod_i r_i = 0\right\}$ goes to zero as $\rho_2$ goes to infinity. 
\qed
\begin{teo}
\label{NGH2}
For any $\epsilon > 0$ there exist $R$ such that for all $\rho_2 > R$ and for  $\rho_1\rho_2$ small enough:\\
1) The fibres of the map $\pi_2$ have diameter at most $\epsilon$ .\\
2)  The Riemannian manifold $\mathbb{X}^{\n -1}_{k_1,k_2}$ is at a (standard) Gromov-Hausdorff  distance of less than $\epsilon$  from its image under the map $\pi_2$\\
3) The  image of $\mathbb{X}^{\n -1}_{k_1,k_2}$ under the map $\pi_2$ is at a distance less than $\epsilon$ from the hypersurface $\prod_{i=0}^\n z_i = 0$ inside $H^\n$, with respect to the distance induced by the Fubini-Study metric.\\
4) The  image of $\mathbb{X}^{\n -1}_{k_1,k_2}$ under the map $\pi_2$ is at a normalized distance less than $\epsilon$ from the hypersurface $\prod_{i=0}^\n z_i = 0$ inside $H^\n_{\rho_1}$, with respect to the (degenerate) distance on $H^\n$ induced by the projection map from $(\cnum^*)^\d$  and the metric on this  space induced by the map $\phi_{k_1,k_2}$ as described in Lemma ~\ref{defphi}
\end{teo}

\textit{Proof}
As in the previous theorem, we will not seek optimal constants for the approximation, and we will content ourselves with proving the statements for the space $\tilde{\mathbb{X}}^{\n -1}_{k_1,k_2}$ and the map $\tilde{\pi}_2$ instead. This will clearly prove also  the analogous statements for $\mathbb{X}^{\n -1}_{k_1,k_2}$ and $\pi_2$. 

1) This is clear as the diameter of the fibres of $\tilde{\pi_2}$ is bounded by a constant multiple of $\rho_1$.

2) This follows immediately from the previous point.

3) 
Define for $\lambda_1,\delta\in\mathbb{R}^+$ and $\lambda_2 > \d$ the compact Riemannian manifolds
\[U^{\delta}_j =  
\left\{[\bar{r},\bar{\eta}] \in 
\tilde{\pi}_2\left(\tilde{\mathbb{X}}^{\n -1}_{k_1,k_2}\right)~|~ ~
\forall i\not= j~r_i  > \delta \right\}, ~ ~
U^\delta = \bigcup_j U^{\delta}_j\]
\[V^{\delta}_j =  
\left\{[\bar{r},\bar{\eta}] \in 
H^\n~|~ r_j = 0,~
\forall i\not= j~r_i  > \delta \right\}, ~ ~
V^\delta = \bigcup_j V^{\delta}_j\]
Notice that for a point in $U^\delta_j$ you have 
\[e^{-4\pi^2\rho_2^2 r_j^2} = 1 - \sum_{i\not= j}e^{-4\pi^2\rho_2^2 r_i^2} > 1 - \frac{\n}{e^{4\pi^2\rho_2^2 \delta^2}}\]
and therefore for $\delta$ small enough $r_j$ can be forced to be an arbitrarily small fraction of one. Notice also that increasing $\rho_2$ only forces $r_j$ to be even smaller. This shows that any point in $U^\delta$ is at a distance from some point in  $V^\delta$ which can be made to be arbitrarily small by taking $\delta$ to zero.
Now pick any point 
$[\bar{r},\bar{\eta}] \in V^{\delta}_j$. For $\delta$ small enough, we can find a point $(\tilde{r}_0,...,\tilde{r}_\n,\bar{\eta})$ such that 
$\tilde{r}_j = 0$, 
$\forall i\not= j~ \tilde{r}_j > \delta $, 
$\sum_i(\tilde{r}_i)^2 =(1-\delta)^2 $ and for all $i$
$(\tilde{r}_i - r_i)^2$ is a small fraction of one. This can be done by putting $r_j$ to zero, and then rescaling the largest $r_i$ to adjust the value of the radius. For $\rho_2$ large enough, we then have that 
\[(1-\delta)^2 <\sum_{i\not= j}(\tilde{r}_i)^2 -\frac{1}{4\pi^2}log\left(1 - \sum_{i\not= j}e^{-4\pi^2\rho_2^2(\tilde{r}_i)^2}\right) < 1\]
Finally, we conclude that there must be an $s < \delta $ such that
\[ \sum_{i\not= j}(\tilde{r}_i + s)^2 -\frac{1}{4\pi^2}log\left(1 - \sum_{i\not= j}e^{-4\pi^2\rho_2^2(\tilde{r}_i+s)^2}\right) = 1\]
and therefore the point 
\[[\tilde{r}_0 + s,...,\tilde{r}_{j-1} + s,\sqrt{ -\frac{1}{4\pi^2}log\left(1 - \sum_{i\not= j}e^{-4\pi^2\rho_2^2(\tilde{r}_i+s)^2}\right)},\tilde{r}_{j+1} + s,..., \tilde{r}_\n + s,\bar{\eta})]\]
lies in $\mathbb{X}^{\n -1}_{\lambda_1,\lambda_2}$ and has a  distance less than  $8(\d)^3\delta$ from the original point $[\bar{\theta},\bar{r},\bar{\eta}] \in V^{\delta}_{j}$. \\
The  argument above shows that for $\delta$  small enough, $\rho_1\rho_2$ small enough and $\rho_2$ large enough the distance 
$d(V^{\delta}, \mathbb{X}^{\n -1}_{k_1,k_2}) $ can be made as small as we want.\\
4) 
The metric is  induced on the quotient by the metric
\[ \sum_{i=0}^\n\frac{16\pi^4\rho_1^2\rho_2^2r_i^2}{e^{4\pi^2\rho_2^2r_i^2}}(dr_i)^2 + 
\frac{e^{4\pi^2\rho_2^2r_i^2}}{4\pi^2 \rho_2^2\rho_2^2}(d\eta_i)^2\]
Then it is clear that for $\rho_1\rho_2$ small enough and $\rho_2$ large enough, the "vertical" directions (those expressed in terms of the $\frac{\partial}{\partial \eta_i}$) become dominant, and therefore from the proof of the previous point we see that as the movements needed were all orthogonal to these directions, the proof follows.
 \qed
 \begin{rmk}
The complex structure on $H^\n$ associated to the metric described in point 4 of the theorem and $\rho_2^2$ times the Fubini-Study two-form is $J_{\rho_1,\rho_2}$, described in Definition ~\ref{dfn:cpxhn} and Proposition ~\ref{prokaonhn}. This structure has a "large complex structure limit" for $\rho_1,\rho_2$ as in the statement of the theorem.
\end{rmk} 
 \begin{pro}
 \label{antic-smooth-2}
 The hypersurface $\prod_i z_i = 0$ inside $H^\n$ can be approximated arbitrarily well by a Calabi-Yau hypersurface inside $H^\n_{\rho_2}$, with respect to the Fubini-Study metric (and also with respect to the metric described in point $4$ of the preceding theorem)
 \end{pro}
 The proof for the Fubini-Study metric is an easy variant of the proof of Proposition ~\ref{antic-smooth-1}, while the proof for the second metric is an easy variant of the proof of point $4$ of the previous theorem. 
 \qed 

\begin{teo}
1)
For $e^{4\pi^2\rho_2^2}\to (\d)^-$ the limiting manifolds are of the form
\[\left(\Torus^\d/N_{\Delta_\n}\right) \times \left(\Torus^\d/(N_{\Delta_\n^*})\right)\]
with flat metrics. The flow $\alpha_t$ induces a flow on the boudary component $e^{4\pi^2\rho_2^2} =  (\d)$. Moreover, the metric and the form $\omega_D$ induce a \ka structure on the limiting manifold.\\
2) For $\rho_2$ fixed, and $\rho_1\to +\infty$ the limiting manifolds are  torus fibrations over $S^{\n -1}$ with respect to normalized distance. Their diameter diverges, and the metric induced on $S^\n$ is  different from (a constant multiple of) the standard one.
3) For $\rho_2$ fixed, and $\rho_1\to +0$ the limiting manifolds are flat tori $\Torus^\n$
\end{teo}
\textit{Proof}
1) This follows from the fact that for $e^{4\pi^2\rho_2^2}\to (\d)^-$ the set of $\bar{r}$ which satisfy both $\mu_1(\bar{r}) = k_1$ and $\mu_2(\bar{r}) = k_2$ reduces to a single point. the statement on the action of $\alpha_t$ and on the \ka structure are clear from the respective definitions.
2) The only observation to make is that for $\lambda_1>>0 ,\frac{1}{\lambda_2} >> 0$ the smooth manifold
\[\left\{(r_0,...,r_\n)\in(\rnum^+)^\d~|~\sum_ir_i^2 = \lambda_1^2,~\prod_i\lambda_1 = \lambda_2\right\}\]
is diffeomorphic to $S^{\n -1}$. One way to verify this is to observe that it retracts to
\[\left\{(r_0,...,r_\n)\in(\rnum^+)^\d~|~\sum_ir_i^2 = \lambda_1^2,~\sum_i\rho_i = \frac{\lambda_1}{2}\right\}\]
along the geodesics originating from $\frac{\rho_1}{\sqrt{\d}}(1,....,1)$. The rest is clear from the description of the action of $\alpha_t$ on the parameters $\rho_1,\rho_2$.\\
3) This is clear from the description of the action of $\alpha_t$ on the parameters $\rho_1,\rho_2$.
\qed
\begin{rmk}
For $\rho_2\to \infty$ with $\rho_1$ fixed it seems that the objects that one obtains are too wildly singular: notice that if a single $r_j\to 0$ the fibre $\omega_1^0$ diverges to a large $\Torus^1$, and the other directions remain with limited diameter. If however two radii go to zero, the fibre diverges to a large $\Torus^2$, and so on, with a diverging behaviour which depends on the "relative speeds of approach" to zero of the various radii.
\end{rmk} 
\section{The general hypersurface case and conclusions}
In the present paper we tried to remain focused on the final aim of building a family of geometric objects interpolating between a family of Calabi-Yau manifolds and its mirror. We obtained a family which interpolates the large \ka structure limit point in the Calabi-Yau moduli space with {\em some} large complex structure limit point. It is not clear if this limit point is the one  expected on physical grounds (and described, for instance, in ~\cite{CDGP}). Even if the answer to this question is negative however, it may simply be that we concentrated on the wrong part of the deformation space of the WSD manifolds $\mathbb{X}^m$. Indeed, the deformation space that we obtained, by varying $\rho_1$ and $\rho_2$ (or equivalently $k_1,k_2$) is but a small part of a much larger deformation space. We believe that the deformation space of $\mathbb{X}^m$ should at least contain a nondegenerate WSD manifold of dimension $3m$. One way to obtain some of these extra deformations is to perform a transformation similar to $\alpha_t$, but this time varying by different parameters $t_0,...,t_\n$ in the various directions associated to the decomposition of $(\cnum^*)^\d$. Another source of deformations are the transformations $\beta_t$ (described in ~\cite{G2}), which can be generalized in the same way as the $\alpha_t$ to provide a $\d$-dimensional family of parameters. A third source of deformations comes from a twisting of the $\Torus^\n$ fibrations. We think that all these deformations are not independent, but generate enough independent deformations to span the WSD manifold mentioned above. Actually it may be that the WSD manifold of deformation has also $2$ extra degenerate dimensions associated to the degenerate directions of $\mathbb{X}^m$, which would bring its total dimension to equal that of $\mathbb{X}^m$. It will be interesting to study the geometry of the whole deformation space, and to see how it relates to the geometry of the deformation spaces of the Calabi-Yau manifolds and their mirrors.

The construction of the manifold $\mathbf{X}^m_{k_1,k_2}$ can be  generalized to include the case of polytopes which have a property that we describe in the following. We use the same notations $\Delta,\Delta^*,F_\Delta, F_{\Delta^*}, f_\Delta, f_{\Delta^*}, D_{\Delta}, D_{\Delta^*}$ as in section ~\ref{sec:prelim}. We indicate with $n$ the dimension of the span of $\Delta$, and with $d$ the number of its vertices.
\begin{dfn}
The polytope $\Delta$ has property SD if the following hold:\\
1) $\Delta$ is integral, and so is its dual polytope $\Delta^*$.\\
2) The number of vertices and the dimension of the spanned space for $\Delta$ and for $\Delta^*$ are the same\\
3) The subgroup $Ker\left(f_{\Delta}f_{\Delta^*}^*\right) = Ker\left(f_{\Delta^*}f_{\Delta}^*\right)$ of $\Torus^n$ is finite
\end{dfn}
Once the above property holds, the construction that we did to build the $\mathbb{X}^\n$ carries over, and we obtain a $2(d-n)$ dimensional family of WSD manifolds of dimension $ d - 2(d-n) + 2n = 4n - d$.

The following remarks are addressed to physicists. Mathematicians can safely skip them.
The first question which arises from a physicist's perspective is how do we recover the $B$-field from this approach.
The manifolds $\mathbb{X}^m$ are probably associated to vanishing $B$-field on the limiting manifolds. One possible way to obtain a nonzero $B$-field should be to start the construction with a different section $\tilde{\sigma}$: not the flat one that we used, but one which acquires monodromy as you go around the homology cycles of the basis. At least this is the way in which the $B$-filed shows up when doing the construction for elliptic curves (cf. ~\cite{G2}). Therefore there should be $\n - 1$ deformation directions associated with varying the $B$-field. One first attempt to verify this would be to perform the computations in the case $\n =1$, where we have a reasonably clear picture of the origin of the $B$-filed from ~\cite{G2}.

The second question which seems relevant to physics is wether the "numerology" of dimensions has some meaning, or is just a coincidence: recall that in the (physically most relevant) case of $\n = 4$ we end up with an $11$-dimensional Riemannian manifold. One possible way to find some meaning in these number would be to build a gauge theory on the bundle $\wedge^*T^*\mathbb{X}^3$ of all exterior powers of the cotangent bundle of $\mathbb{X}^3$, using the generalized lagrangian dynamics introduced in ~\cite{G1}. That part of that paper was precisely aimed at building a generalization of lagrangian (and hamiltonian) dynamics well suited to study PDE's. It may be a coincidence, but the object which was the  outcome of the theory in the case of PDE's with two independent variables was {\em exactly} (a piece of) a WSD structure. Another intriguing fact is  that this hypothetical theory would be a $11$ dimensional gauge theory which on the boundary of its deformation space would give rise to the same theory which comes out of a $\sigma$-model. It is very tempting to conjecture that if it actually exists this theory is  (strongly related to) what is usually called $M$-theory. In this respect one should try to build a representation of $E_8$ on the above mentioned bundle, and again there is some interesting "numerology" coming out of the computation of its rank. We will investigate these issues in a future paper.
\bigskip
~  ~\\
Michele Grassi (grassi@dm.unipi.it)\\
Dipartimento di Matematica, Universit\`{a} di Pisa\\
Via Buonarroti, 2\\ 
56100 Pisa, Italy

\end{document}